\documentclass[12pt]{article}
\usepackage[affil-it]{authblk}

\usepackage[absolute,overlay]{textpos}

\usepackage{graphicx}
\usepackage{ragged2e}
\usepackage{tabto}

\usepackage{amsmath}
\usepackage{amsthm}
\usepackage{amsfonts}

\usepackage{float}
\usepackage{multirow}
\usepackage{tikz}
\usepackage{graphics}
\usepackage{color}
\usepackage{url}
\usepackage{afterpage}
\usepackage{multicol}
\usepackage[font=footnotesize]{caption}
\usepackage[font=footnotesize]{subcaption}
\usepackage{marginnote}
\usepackage{tcolorbox}

\usepackage{lmodern}
\usepackage{pgfplots}

\usepackage{bm}
\usepackage{booktabs}
\usepackage{siunitx}
\usepackage{isomath}

\usepackage{url}

\usepackage[plainpages=false, colorlinks=true, citecolor=blue, filecolor=blue,
   linkcolor=blue, urlcolor=blue]{hyperref}
\usepackage{enumitem}

\usepackage{comment}

\usepackage{multicol}
\newcommand{\lp}{\left(}
\newcommand{\rp}{\right)}

\newcommand{\intt}{{\operatorname{int}}}
\newcommand{\extt}{{\operatorname{ext}}}
\newcommand{\nau}{^{(0)}}
\newcommand{\one}{^{(1)}}

\newcommand{\Tint}{T_{\intt}}
\newcommand{\Text}{T_{\extt}}

\newcommand{\rint}{\rho_{\intt}}
\newcommand{\rext}{\rho_{\extt}}

\newcommand{\pint}{p_{\intt}}
\newcommand{\pext}{p_{\extt}}

\newcommand{\partialz}{\frac{\partial}{\partial z}}
\newcommand{\partialzp}{\frac{\partial}{\partial z'}}

\newcommand{\res}{{\text{res}}}

\newcommand{\dz}{{\,\rm{d}}z }

\begin{document}

\baselineskip=15pt

\title{A Model of a Buoyancy-Driven Heat Exchanger, with Implications for Optimal Design}

\author{Sylvie Bronsard \thanks{Courant Institute of Mathematical Sciences, New York University,
New York, NY 10012, 
\texttt{sab980@nyu.edu}} \ and Charles S. Peskin \thanks{Courant Institute of Mathematical Sciences, New York University,
New York, NY 10012,
\texttt{csp1@nyu.edu}}}

\maketitle

\section*{Abstract}

In this paper, we introduce a model for a buoyancy-driven, air-to-air heat exchanger. This model, derived from first principles, features a conservative boundary condition at inflow based on the compressible Bernoulli equation, and a dissipative boundary condition at outflow based on pressure continuity. We solve for the steady-state behavior numerically and asymptotically, with excellent agreement between the two, and we study the tradeoff between the efficiency and air flow predicted by the model.

\section{Introduction}	
It takes significant energy expenditure to maintain comfortable indoor conditions in a well ventilated building \cite{MonahanPowell2011}. In cold climates, fresh air from outdoors must be heated to indoor temperature; the energy required to do this can be lessened by using a device, called a heat exchanger, which uses the exhaust air to warm the incoming air.

Heat exchangers are used both in industry and in nature. For example, they improve fuel efficiency in
airplanes by using heat from the engine to preheat the fuel. They are also used in the making
of wine and beer, for pasteurizing dairy products, and more generally in the pharmaceutical,
food, and drink industries \cite{Taler2019}. In nature, for example, flamingoes and other wading birds use
heat exchangers to keep their bodies warm while standing in cool water.

In buildings, heat exchangers can be active, using a fan or pump to force fluid flow; or passive. Passive heat exchangers make use of the fact that warm air is lighter than cold air. A simple example is a device consisting of a tube with a partition down the middle, allowing warm air to rise through one side of the heat exchanger, and cold air to
sink down the other. As the warm air rises, some of its heat is transferred to the incoming colder air across the partition, thus decreasing heating costs. Although active heat exchangers are more common, passive exchangers do
 not require input energy (so they continue to run during a blackout),
 are more quiet, require less maintenance, and have a longer operating
 life than active heat exchangers \cite{OCH16}.

Buoyancy driven ventilation has been installed in schools \cite{LIPINSKI2020}, with the additional benefit of being particularly effective at expelling small droplets and airborne particles from the indoor air, thus circulating fresh air in such a way as to reduce the spread of infectious disease \cite{bhagat2020}. Since a human being's internal temperature is warmer than the air temperature, the air we exhale naturally rises, carrying with it the particles and pathogens contained within. Buoyancy driven ventilation takes advantage of this natural stratification, and replaces the stale air without forced mixing.

Since heat exchangers are in use in many industrial applications, they have naturally been extensively studied by engineers. This existing literature may be broadly divided into two categories: the first and largest category contains articles and review papers that use experimental methods and simple models based on empirical laws. In these models, the pressure difference determines the fluid velocity, which in turn determines the heat transfer between the fluids, without differential equations describing the fluid flow within the device. These works ask interesting empirical questions about the effect of turbulent versus laminar flow on the fluid velocity and heat transfer between fluids. They study the effects of modifying the shape of the heat exchanger (for example, by putting nails through the partition), and build experimental setups to test their models \cite{Schultz93, AHS11, LS13}. Most, but not all, of the existing literature focuses on driven heat exchangers rather than passive ones.

There has also been recent work using detailed three-dimensional
computational fluid dynamic (CFD) models to simulate fluid flow in driven, aluminum plate-fin heat exchangers, studying fully developed flow conditions and taking into account entrance effects on the fluid flow \cite{Grespan2025}. 

The present paper introduces a framework that is intermediate in
 complexity between the two types of models described above.  Our
 framework is based on one-dimensional compressible steady-state fluid dynamics within a pair of tubes, with heat conduction across a
 partition that separates the two tubes.  We neglect fluid viscosity
 and heat conduction within each tube.  At the entrance to each tube we use the compressible Bernoulli equation, which is energy
 conserving.  At the exit of each tube, we use a continuous-pressure
 boundary condition that dissipates kinetic energy, and this dissipation is essential for the existence of a steady state in our
 model (see Appendix).  Our model is passive, with flow driven by gravity as a result of the difference in density between warm air and
 cold air.  The equations of the model are solved
 numerically and an asymptotic analysis is also done, with excellent agreement between the numerical and asymptotic results.  By both
 methods, we determine the efficiency of the heat exchanger, suitably
 defined, and the relative mass flux, also suitably defined, of fresh air brought in by the heat exchanger.  We observe that the objective of maximizing efficiency is in conflict with the objective
 of maximizing the relative mass flux, but we are able to maximize the minimum of these two objective functions and in that sense to
 optimize the design of a passive heat exchanger.

\section{Mathematical Formulation} \label{section2}

We model a buoyancy-driven heat exchanger as a pair of one-dimensional vertical tubes in contact along their length, allowing warm air to rise through one side of the heat exchanger, and cold air to sink down the other. We ignore the effects of viscosity and thermal conductivity within the fluid, so that heat is only transferred across the partition between the two tubes. We focus on the steady-state behavior of the heat exchanger.

In other words, we want to study steady gas dynamics in a vertical tube, under the effect of gravity, and with a heat source or sink along the length of the tube. The partition between the two tubes is assumed to be a rectangle of width $W$ and height $H$.
%later: although we picture it like this, it would be equivalent to concentric annuli with circumpherence W and height H

A diatomic gas at ordinary temperatures satisfies the ideal gas law, 

\begin{equation}
p = \rho\lp\frac{RT}{m}\rp, \label{idealgaslaw}
\end{equation}
and has internal energy (per mole)
\begin{equation}
e = \frac{5}{2}\lp\frac{RT}{m}\rp, \label{internalenergy}
\end{equation}
where $p$ is the pressure, $\rho$ the density, $T$ the temperature, $m$ the mass per mole, and where R is the molar gas constant.

\begin{figure}
\centering
\includegraphics[width=6cm]{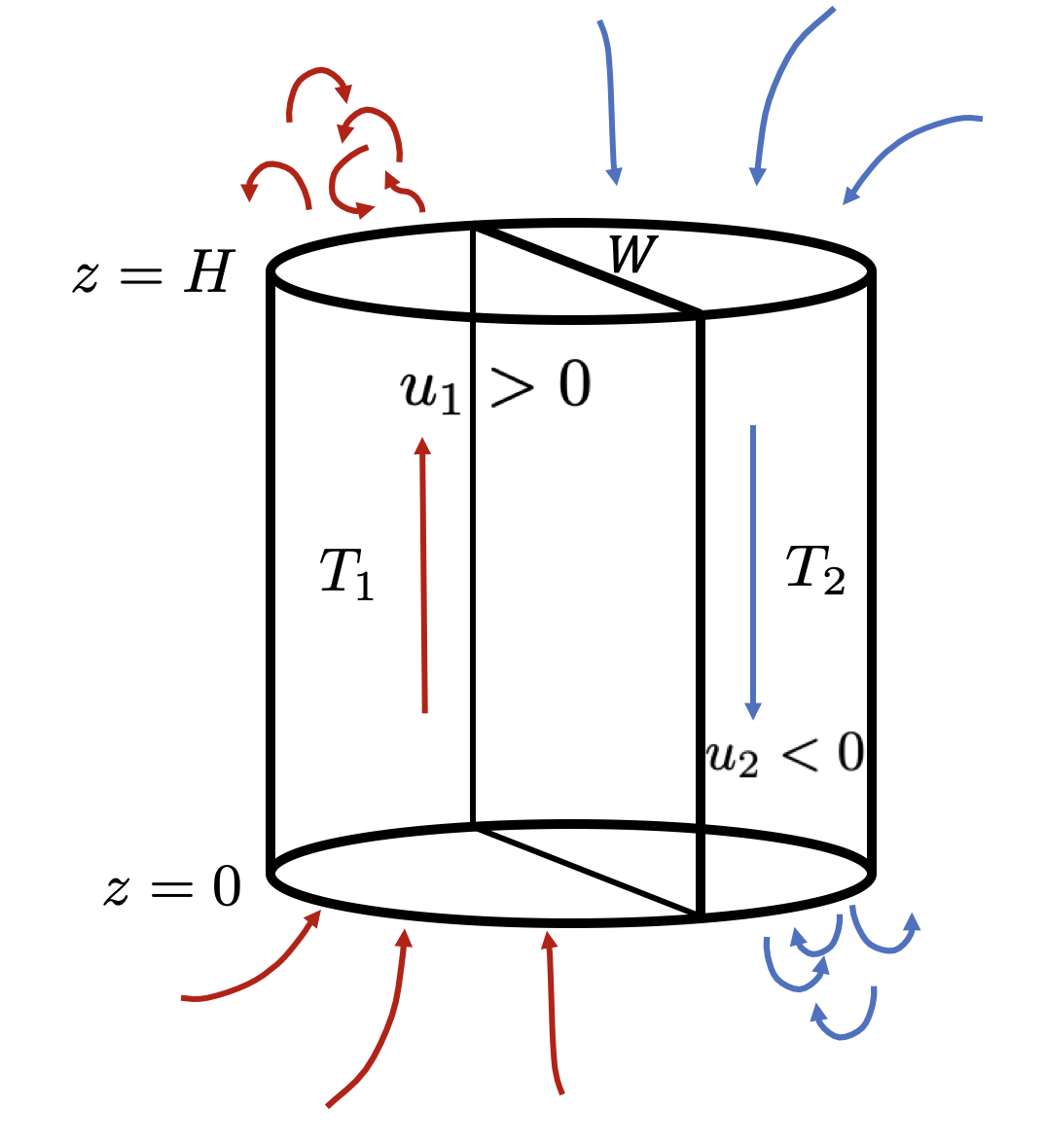}
\caption{Warm air rises through the left side and sinks down the right side, with heat transfer across the rectangular partition. We assume a smooth, energy-preserving flow at the inflow, and a turbulent, energy-dissipating flow at the outflow.}
\label{fig:setup}
\end{figure}

We use the subscript $i=1,2$ to denote the two tubes, and choose tube 1 to be the tube through which warm air rises and tube 2 as the tube through which cool air descends.\footnote{Since the flow is buoyancy driven, not forced, the heat exchanger decides for itself which tube will be the outflow and which the inflow. We merely label the outflow tube as tube 1.} So $T_i, u_i, \rho_i$, and $p_i$ are, respectively, the temperature, velocity, density, and pressure in tube $i$. Our spatial coordinate $z$ denotes height measured from the bottom of the heat exchanger, and our definition of the (steady) fluid velocity $u_i(z)$ in each tube is such that $dz/dt = u_i(z(t))$ for a fluid particle
 with trajectory $z(t)$ in tube $i$.  Thus $u_1(z) > 0$, and $u_2(z) < 0$. We assume for simplicity that both tubes have the same constant cross-sectional area $A$.  See Figure \ref{fig:setup} for an illustration of the setup as described.

The equations of mass, momentum, and energy conservation in each tube are then formulated as follows.
\begin{align}
\partialz \lp \rho_i u_i A\rp=0, \label{masscons}\\
\partialz \lp p_i A + u_i A \rho_i u_i \rp &= -\rho_i g A, \label{momcons}\\
\partialz \lp p_i A u_i + ( e_i \rho_i + \frac12 \rho_i u_i^2) A u_i\rp &= q_i W-\rho_i g A u_i,
\label{encons}
\end{align}
where $i=1$ or 2, where $q_2 = -q_1 = \sigma (T_1-T_2)$, where $\sigma$ is the thermal conductivity divided by the thickness of the partition, and where $W$ is the width (area per unit height) of the partition. These equations hold for $z \in [0,H]$, where in both tubes $z = 0$ is in the house and $z = H$ is exposed to the outside air.

We can rewrite these equations in a simpler form. First, we introduce the mass flux per unit cross-sectional area, $\Phi_i = \rho_i u_i$, and assume that the air flow through the walls of the house is negligible compared to that through the heat exchanger, so $\Phi_1 = -\Phi_2 =: \Phi > 0$. The conservation of mass equation \eqref{masscons} then reduces to the statement that $\Phi$ is constant.

We also note that 
\begin{equation}
\frac{q_i}{\rho_iu_i}= -\frac{q}{\Phi},
\end{equation}
where $q=q_2 = -q_1$.

Using this, we divide our remaining conservation equations by $A$ and rewrite them in terms of $\Phi$ instead of $u_i$. We then have:
\begin{align}
\partialz \lp p_i + \frac{\Phi^2}{\rho_i} \rp &= -\rho_i g, \label{momentumcons}\\
\partialz \lp \frac{p_i}{\rho_i} + e_i + \frac12 \lp\frac{\Phi}{\rho_i}\rp^2\rp &= - \frac{q}{\Phi}\frac{W}{A}-g,
\end{align}
We can further simplify the energy conservation equation. From equations \eqref{idealgaslaw} and \eqref{internalenergy},
\begin{equation}
\frac{p_i}{\rho_i} + e_i  = \frac72 \frac{RT_i}{m},
\end{equation}
and we rewrite $u_i^2$ as $\frac{\Phi^2}{\rho_i^2}$, so that we are left with the following conservation of energy equation:

\begin{equation}
\partialz \lp \frac72 \frac{RT_i}{m} + \frac12 \frac{\Phi^2}{\rho_i^2}\rp = - \frac{q}{\Phi}\frac{W}{A}-g. \label{energycons}
\end{equation}

\subsection{Boundary Conditions} \label{bc}
We now turn to the boundary conditions. At the inflow to each tube, we picture steady flow from a reservoir; the interior of the house for tube 1, the exterior air for tube 2. We assume adiabatic equilibrium between the air at the tube entrance and the air in the reservoir, which we take to be at the same height as the tube entrance. We assume that the air flow from the reservoir to the entrance of the corresponding tube is governed by the 3D compressible Euler equation,

\begin{equation}
\rho \lp \frac{\partial u}{\partial t}+u\cdot \nabla u \rp+\nabla p =0. \label{gen_mom_cons}
\end{equation}
%For a parcel of gas in which there is no heat flux across the boundary of the parcel, we also know that the adiabatic equation of state holds:
with the adiabatic equation of state
 \footnote{since we are assuming that there is negligible heat transfer
 (as well as negligible viscosity) between adjacent parcels of fluid
 within the flow from each reservoir into the corresponding tube.}

\begin{equation}
p/\rho^{\gamma}= \text{constant,} \label{eqnofstate}
\end{equation}
where $\gamma=\frac75$ is the ratio of the specific heat at constant pressure to the specific heat at constant volume for a diatomic gas. 

Since all parcels of fluid come from the same place, namely the reservoir, the constant in equation \eqref{eqnofstate} is independent of position throughout either one of the two entrance flows. The constants may be different, however, for the two entrance flows. In particular, we can write

\begin{equation}
p/\rho^{\gamma}= p_{\text{res}}/\rho_{\text{res}}^{\gamma}, \label{adi_eqn_state}
\end{equation}
where $p_\res$ and $\rho_\res$ are, respectively, the pressure and density in the reservoir.

Armed with equations \eqref{gen_mom_cons} and \eqref{adi_eqn_state}, we are ready for action. We first divide \eqref{gen_mom_cons} by $\rho$ to obtain

\begin{equation}
\frac{\partial{u}}{\partial{t}}+u\cdot \nabla u +\frac{\nabla p}{\rho} =0.
\end{equation}
We can rewrite $\frac{\nabla p}{\rho}$ as follows:
\begin{align}
\nonumber \frac{\nabla p}{\rho} &= \frac{\rho_\res}{\rho_\res}\cdot \frac{p_\res}{p_\res}\cdot \frac{\nabla p}{\rho} \\
&= \frac{p_\res}{\rho_\res}\cdot \frac{\rho_\res}{\rho} \cdot \nabla\lp\frac{p}{p_\res}\rp.
\end{align}
Here we note that $\frac{p_\res}{\rho_\res} = \frac{RT_\res}{m}$, and that \eqref{adi_eqn_state} is equivalent to $\frac{\rho_\res}{\rho}=\lp\frac{p_\res}{p}\rp^{1/\gamma}$, so
\begin{align}
\nonumber \frac{\nabla p}{\rho} &=\frac{RT_\res}{m}\lp\frac{p}{p_\res}\rp^{-1/\gamma}\nabla\lp\frac{p}{p_\res}\rp\\
&=\frac{RT_\res}{m}\lp\frac{\gamma}{\gamma-1}\nabla\lp\lp\frac{p}{p_\res}\rp^{\frac{\gamma-1}{\gamma}}\rp\rp.
\end{align}
Using \eqref{adi_eqn_state} once again, we further note:
\begin{equation}
\frac{p}{\rho}= \frac{p}{\rho_\res\lp\frac{p}{p_\res}\rp^{\frac{1}{\gamma}}} =\frac{p_\res}{\rho_\res}\lp\frac{p}{p_\res}\rp^{\frac{\gamma-1}{\gamma}}=\frac{RT_\res}{m}\lp\frac{p}{p_\res}\rp^{\frac{\gamma-1}{\gamma}},
\end{equation}
so $\frac{\nabla p}{\rho} = \frac{\gamma}{\gamma -1}\nabla\lp\frac{p}{\rho}\rp$. We can therefore write

\begin{equation}
\frac{\partial{u}}{\partial{t}}+u\cdot \nabla u +\frac{\gamma}{\gamma -1}\nabla\lp\frac{p}{\rho}\rp =0.
\end{equation}
In steady flow $\frac{\partial{u}}{\partial{t}}=0$, and for a diatomic ideal gas $\gamma = \frac75$, therefore
\begin{align}
&u\cdot \nabla u +\frac72\nabla\lp\frac{p}{\rho}\rp =0\\
\implies &u\cdot \lp u\cdot\nabla u +\frac72\nabla\lp\frac{p}{\rho}\rp\rp =0\\
\implies &u\cdot\nabla \lp \frac12 \|u\|^2 +\frac72\lp\frac{p}{\rho}\rp\rp =0,
\end{align}
since $u\cdot \lp u\cdot\nabla u\rp =u\cdot\nabla \frac12 \|u\|^2$.
We conclude that the quantity $\frac12 \|u\|^2 +\frac72\lp\frac{p}{\rho}\rp$ is constant along streamlines. This is known as the compressible Bernoulli condition \cite[\S 3.11]{clancy1975aerodynamics}. Since $u=0$ in the reservoir, at each tube entrance we then have $\frac12 \|u\|^2 +\frac72\lp\frac{p}{\rho}\rp = \frac72\lp\frac{p_\res}{\rho_\res}\rp$. The entrance to tube 1 is at $z = 0$, and the entrance to tube 2 is at $z=H$. Therefore,
\begin{align}
\frac12 u_1^2(0) + \frac72 \frac{RT_1(0)}{m}  &= \frac72 \frac{R\Tint}{m}\label{bc1}\\ 
 \frac12 u_2^2(H) + \frac72 \frac{RT_2(H)}{m} &= \frac72 \frac{RT_{\extt}}{m} \label{bc2}
\end{align}
and
\begin{align}
\frac{p_1(0)}{\rho_1(0)^\frac57}&=\frac{\pint}{\rho_{\intt}^{\frac57}} \\
\frac{p_2(H)}{\rho_2(H)^\frac57}&=\frac{p_{\extt}}{\rho_{\extt}^{\frac57}}.
\end{align}
Although we expect a smooth, energy-preserving flow at the entrance of each tube, as the fluid leaves the tube it collides with the still air outside the tube, and so its kinetic energy is lost. As a result, there is not the pressure recovery that would be predicted by the Bernoulli equation; instead, the pressure in the tube at the exit matches the pressure of the external air: 
\begin{align}
&p_1(H) = p_0 \label{outbc1}\\
&p_2(0) = \pint. \label{outbc2}
\end{align}
Here, $p_0$ is the atmospheric pressure and $\pint$ is the pressure in the house. Although we regard $p_0$ as given, $\pint$ is an unknown which is determined by the condition of mass conservation. Since we are assuming no flow through the walls of the house, the net flow of mass through the heat exchanger as a whole must be zero, and $\pint$ will have to adjust to achieve this condition.

Equations $\eqref{outbc1}-\eqref{outbc2}$ imply energy dissipation at the outflow of each tube, since kinetic energy is dissipated. This will be shown in the appendix, and it is an important feature of our model, needed to achieve a steady state since there is no other dissipative mechanism in the model.

\subsection{Power and Efficiency} \label{eff}

We are interested in the energy efficiency as well as the mass flux of fresh air predicted by the above model. In this section we define the power $P$ required to maintain the internal temperature of the house at $\Tint$:

\begin{equation}
P = \sum_{i=1}^2 \lp Au_ip_i + Au_i \lp e_i\rho_i+\frac12\rho_iu_i^2\rp\rp.
\end{equation}
Note that $P$ is independent of $z$.  This follows by summing over $i=1,2$ in equation \eqref{encons} and noting that the sum of the right-hand sides is zero, since $\rho_i u_i$ has the same magnitude and opposite sign in
 the two tubes.  Thus $P$ is the flux of energy from below across any
 plane $z =$ constant.  The term in $P$ involving the pressure is the net rate at which work is being done by the pressure from below the plane
 $z =$ constant on the air above it, and the remaining terms in $P$
 describe the net rate of upward transport by convection, across the
 plane $z =$ constant, of the internal energy and the kinetic energy of
 the air.

The following formula relates the power requirement $P$ to the important
 overall variables: $\Tint$, $\Text$, $\Phi$, and $Q = \frac{W}{A}\int_0^Hq\dz$.

We claim that
\begin{equation}
P = A\Phi\lp\frac72\frac{R}{m}\lp\Tint-\Text\rp-gH-\frac{Q}{\Phi}\rp. \label{pwr}
\end{equation}
%\end{lemma}
\begin{proof}
Since $\Phi_i = \rho_iu_i$ and $\Phi_1 = -\Phi_2 = \Phi$, we can write
\begin{align}
\nonumber P &= \sum_{i=1}^2 A\Phi_i \lp \frac72\frac{RT_i}{m}+\frac12u_i^2 \rp \\
&= A\Phi  \lp  \frac72 \frac{R}{m}\lp T_1-T_2\rp +\frac12\lp u_1^2-u_2^2\rp\rp. \label{pow_def}
\end{align}
But we also know, by integrating \eqref{energycons} from $z = 0$ to $z=H$, that 
\begin{align}
\nonumber \left.\frac72\frac{RT_i}{m}\right|_0^H+\left.\frac12\frac{\Phi^2}{\rho_i^2}\right|_0^H &= \int_0^H -\frac{q}{\Phi}\frac{W}{A}-g \dz \\
&= -\frac{W}{\Phi A}\int_0^Hq\dz - gH.
\end{align}
Using our boundary condition \eqref{bc1}, we now have that

\begin{equation}
\frac72\frac{R}{m}T_1(H)+\frac12u_1^2(H) = \frac72 \frac{R\Tint}{m}-gH-\frac{Q}{\Phi}.
\end{equation}
Subtracting from this the boundary condition \eqref{bc2} gives us

\begin{equation}
\frac72\frac{R}{m}T_1(H)+\frac12u_1^2(H) - \lp \frac72 \frac{RT_2(H)}{m}+\frac12u_2^2(H)\rp = \frac72 \frac{R}{m}(\Tint-\Text) -gH - \frac{Q}{\Phi},
\end{equation}
or, written in a more evocative form,

\begin{equation}
\frac72\frac{R}{m}\lp T_1(H)-T_2(H)\rp +\frac12\lp u_1^2(H)-u_2^2(H)\rp = \frac72 \frac{R}{m}(\Tint-\Text) -gH - \frac{Q}{\Phi}.
\end{equation}

We now turn back to our expression for the power requirement of a heat exchanger, given by equation \eqref{pow_def}. Since the right hand side of \eqref{pow_def} does not depend on $z$, we can evaluate it at any height, in particular at $z = H$:
\begin{align}
\nonumber P &= A\Phi  \lp  \frac72 \frac{R}{m}\lp T_1(H)-T_2(H)\rp +\frac12\lp u_1^2(H)-u_2^2(H)\rp\rp\\
&= A\Phi\lp\frac72 \frac{R}{m}(\Tint-\Text) -gH - \frac{Q}{\Phi}\rp,
\end{align}
which is the same as \eqref{pwr}.\\
\end{proof}

The purpose of a heat exchanger is to bring in fresh air in an energy-efficient manner, so we are interested in the energy cost per unit mass of fresh air that is brought in; that is, in the quantity

\begin{equation}
\frac{P}{A\Phi} = \frac72 \frac{R}{m}(\Tint-\Text) -gH - \frac{Q}{\Phi}. \label{pwrperfreshair}
\end{equation}

We can define the efficiency of a heat exchanger as the relative difference in this quantity in comparison to the situation in which there is no heat exchange, i.e. when $Q=0$.
\begin{align}
E :&= \frac{\lp\frac{P}{A\Phi}\rp_{Q=0}-\lp\frac{P}{A\Phi}\rp}{\lp\frac{P}{A\Phi}\rp_{Q=0}}\\
&= \frac{Q/\Phi}{\frac72\frac{R}{m}\lp\Tint-\Text\rp-gH}.\label{effdef}
\end{align}

In equation \eqref{pwrperfreshair} it is reasonable to assume that $P$, $Q$, and $\Phi$ are all positive.  For $P$, this is because the gravity-driven heat exchanger is a passive device, the power requirement of which cannot be negative. For $Q$, this is because heat flows passively from the warmer air (tube 1) to the cooler air (tube 2).  For $\Phi$ it is because the warm air will
 naturally be rising and the cool air descending.  With $P$, $Q$, and $\Phi$ positive, it follows from \eqref{pwrperfreshair} that the denominator in \eqref{effdef} is positive and also that the numerator in \eqref{effdef} is both positive and less than the
 denominator.  Thus, $E\in(0,1)$, as it should be.

 Positivity of the denominator in equation \eqref{effdef} deserves further comment. This is a restriction on the given parameters for our model of the heat exchanger to make sense. What presumably happens when this condition fails is that symmetry breaking does not occur, and the air remains at rest in both tubes. There is then no heat loss, and no fresh air entering
 the house.

\subsection{Dimensionless formulation} \label{dimless}

To make the system easier to work with, we nondimensionalize. We start from the conservation equations in the form \eqref{momentumcons} and \eqref{energycons}. To make the equations dimensionless, we choose units of temperature, length, mass and time. We choose the height of the tube $H$ as the unit length, $\Text$ as the unit temperature, $\pext$ as the unit pressure, and $\rext$ as the unit density. It follows that:
\begin{align}
&\pext H^3= \text{ unit energy;}\\
&\rext H^3 = \text{ unit mass;}\\
&\sqrt{\frac{\pext}{\rext}}=\sqrt{\frac{R\Text}{m}}=\text{ unit velocity};\\
&H/\sqrt{\frac{R\Text}{m}}=\text{ unit time}.
\end{align}

We now introduce dimensionless variables, 
\begin{equation}
T = \Text T', \qquad
p = \pext p', \qquad
\rho = \rext \rho'. \qquad
\end{equation}
Substituting these into the ideal gas law, we have:
\begin{align}
p &= \frac{R}{m}T\cdot \rho \\
\pext p' &= \frac{R}{m}\Text T' \cdot \rext \rho'. 
\end{align}
Since $\pext = \lp R/m\rp \Text\rext$, this can be rewritten as:

\begin{equation}
p' = T' \rho'. \label{idealgas}
\end{equation}
We will also need the following:
\begin{align}
g &= \frac{1}{H} \lp \frac{R\Text}{m}\rp g' \\
\sigma &= \frac{\pext}{\Text}\lp\sqrt{\frac{R\Text}{m}}\rp\sigma' \\
\Phi &= \rext \lp\sqrt{\frac{R\Text}{m}}\rp \Phi',
\end{align}
and, lastly, $z = Hz'$.

Since $z \in (0,H)$, our new variable $z' \in (0,1).$ In terms of our dimensionless variables, equation \eqref{momentumcons} becomes

\begin{equation}
\frac{1}{H}\partialzp \lp \pext p_i' + \frac{1}{\rext\rho_i'}\rext^2 \frac{R\Text}{m}(\Phi')^2\rp = -\rext\rho_i'\frac{1}{H}\frac{R\Text}{m}g'.
\end{equation}
Since $\pext = \rext \frac{R\Text}{m}$, all the dimensional constants cancel out, leaving us with the following momentum conservation equation:

\begin{equation}
\partialzp \lp p_i' + \frac{(\Phi')^2}{\rho_i'}\rp = -\rho_i' g'. \label{adimmom}
\end{equation}

We now turn to energy conservation. Substituting our dimensionless variables into \eqref{energycons} with $q = \sigma(T_1-T_2)$, 

\begin{equation}
\frac{1}{H} \partialzp \lp \frac72 \frac{R\Text}{m}T_i' + \frac12 \frac{1}{\rext^2(\rho_i')^2}\rext^2\frac{R\Text}{m}(\Phi')^2\rp=-\frac{W}{A}\frac{\pext}{\rext\Text}\frac{\sigma'}{\Phi'}\Text(T_1'-T_2')-\frac{1}{H}\frac{R\Text}{m}g'.
\end{equation}
Multiplying the above by $H$ and cancelling $\frac{R\Text}{m}=\frac{\pext}{\rext}$ yields the following.

\begin{equation}
\partialzp \lp \frac72 T_i' + \frac12 \frac{(\Phi')^2}{(\rho_i')^2}\rp = - \lp \frac{WH}{A} \rp \lp\frac{\sigma'}{\Phi'}\rp(T_1'-T_2')-g'.
\end{equation}
Note that $\lp\frac{WH}{A}\rp$ is dimensionless, and this dimensionless parameter only appears as a multiple of $\sigma'$, so it natural to define 

\begin{equation}
\sigma''=\lp\frac{WH}{A}\rp\sigma'.
\end{equation}
This gives us the dimensionless equation,
\begin{equation}
\partialzp \lp \frac72 T_i' + \frac12 \frac{(\Phi')^2}{(\rho_i')^2}\rp = - \lp\frac{\sigma''}{\Phi'}\rp(T_1'-T_2')-g'. \label{a_dim_en}
\end{equation}

We are left with the ODEs \eqref{adimmom} and \eqref{a_dim_en}, which hold for $z' \in (0,1)$ and $i = 1,2$. In these equations, $\Phi'$ is independent of $z'$, and the functions $p_i'(z'), \rho_i'(z'),$ and $T_i'(z')$ are related by the dimensionless ideal gas law \eqref{idealgas}.

We also need to rewrite the boundary conditions in dimensionless form. Substituting and cancelling all dimensional factors, as we have done above with the momentum and energy conservation equations, we obtain:
\begin{align}
\frac72 T_1'(0) +\frac12 \frac{(\Phi')^2}{(\rho_1'(0))^2} &= \frac72\Tint' \\
\frac{p_1'(0)}{(\rho_1'(0))^{\frac75}} &= \frac{\pint'}{(\rint')^{\frac75}} \\
\frac72 T_2'(1) + \frac12 \frac{(\Phi')^2}{\rho_2'(1))^2} &= \frac72 \\
\frac{p_2'(1)}{(\rho_2'(1))^{\frac75}} &=1 \\
p_1'(1)&=1\\
p_2'(0) &=\pint',
\end{align}
where $\Tint' = \Tint/\Text$ is given; and $\pint'=\pint/\pext$ and $\rint'=\rint/\rext$ are unknown but related through $\Tint'$ and the ideal gas law \eqref{idealgas}.

In summary, the behavior of the heat exchanger is governed by the following three dimensionless parameters:
\begin{align}
&\Tint' = \frac{\Tint}{\Text}\\
&g' = \frac{gH}{\lp\frac{R\Text}{m}\rp}\\
&\sigma'' = \lp\frac{WH}{A}\rp \frac{\Text}{\pext\sqrt{\frac{R\Text}{m}}}\sigma.
\end{align}
If we choose reasonable physical parameters (such as setting the indoor air to room temperature, the external air to freezing, the height and width of the heat exchanger to 1m, and so on), the dimensionless parameters will look like $\Tint'=1.0732$, $g'=1.2471\times 10^{-4}$, and $\sigma'' = 0.0018$. We show plots of what the system looks like with those parameters in Section \ref{numexp}.

We also note that the dimensionless parameter $g'$ has an interesting physical interpretation: in the adiabatic atmosphere approximation, the thickness of the atmosphere is given by $\frac{R\Text}{mg}$ \cite{adiabaticatmosphere}. So the dimensionless parameter $g'$ is the ratio of the height of the heat exchanger to the thickness of the atmosphere. 

Lastly, we should also write the efficiency \eqref{effdef} in terms of the dimensionless parameters. Substituting the dimensionless parameters as done above, we find that

\begin{equation}
E' = \frac{Q'/\Phi'}{\frac72\lp\Tint'-1\rp - g'}, 
\end{equation}
where
\begin{equation}
Q' = \int_0^1 \sigma''\lp T_1'-T_2'\rp\dz'.
\end{equation}

From now on, we will always be working with dimensionless parameters, and so will drop the prime notation on the dimensionless variables, but retain them on the parameters. We will also drop one prime on the $\sigma''$ and write $\sigma'$ instead, so that it matches the notation for the other dimensionless parameters.

\section{Numerical Method} \label{section3}
 
In this section, we describe our numerical solver for the above system. Let $\Delta z = \frac{H}{n-1}$, where $n$ is the number of gridpoints we use for our discretization. We discretize our derivatives using a second-order upwind scheme.  Discretizing our conservation of momentum equation gives the following system of equations, in the interiors of tubes 1 and 2 respectively:
\begin{equation}
\frac{3p_{1,j}-4p_{1,j-1}+p_{1,j-2}}{2\Delta z}-\frac{\Phi^2}{\rho_{1,j}^2}\lp\frac{3\rho_{1,j}-4\rho_{1,j-1}+\rho_{1,j-2}}{2\Delta z}\rp - g'\rho_{1,j} = 0, \label{mo_1}
\end{equation}
for $j= 3,...n$, and
\begin{equation}
\frac{-p_{2,j+2}+4p_{2,j+1}-3p_{2,j}}{2\Delta z}-\frac{\Phi^2}{\rho_{2,j}^2}\lp\frac{-\rho_{2,j+2}+4\rho_{2,j+1}-3\rho_{2,j}}{2\Delta z}\rp - g'\rho_{2,j} = 0, \label{mo_2}
\end{equation} 
$j = 1,...n-2$.

For the gridpoint directly downwind of the boundary, we must do something different, as we do not have two upwind points to draw from. We instead use a first order upwind scheme:
\begin{align}
\frac{p_{1,2}-p_{1,1}}{\Delta z}- \frac{\Phi^2}{\rho_{1,2}^2}\lp\frac{\rho_{1,2}-\rho_{1,1}}{\Delta z}\rp - g'\rho_{1,2} = 0 \label{mo_1a}\\ 
\frac{p_{2,n}-p_{2,n-1}}{\Delta z}- \frac{\Phi^2}{\rho_{2,n-1}^2}\lp\frac{\rho_{2,n}-\rho_{2,n-1}}{\Delta z}\rp - g'\rho_{2,n-1} = 0, \label{mo_2a}
\end{align}
This does not affect the overall order of the scheme, which remains second order in the $L^{\infty}$ norm (as shown in the following section, see Fig \ref{fig:convstdy_std}). 

We do the same for our conservation of energy equation:
\begin{equation}
\frac72 \lp\frac{3T_{1,j}-4T_{1,j-1}+T_{1,j-2}}{2\Delta z}\rp -\frac{\Phi^2}{\rho_{1,j}^3}\lp \frac{3\rho_{1,j}-4\rho_{1,j-1}+\rho_{1,j-2}}{2\Delta z}\rp +\frac{\sigma'}{\Phi}\lp T_{1,j}-T_{2,j}\rp +g'=0, \label{en_1}
\end{equation}
for $j= 3,...n$, and
\begin{equation}
\frac72 \lp\frac{-T_{2,j+2}+4T_{2,j+1}-3T_{2,j}}{2\Delta z}\rp -\frac{\Phi^2}{\rho_{2,j}^3}\lp \frac{-\rho_{2,j+2}+4\rho_{2,j+1}-3\rho_{2,j}}{2\Delta z}\rp +\frac{\sigma'}{\Phi}\lp T_{1,j}-T_{2,j}\rp +g'=0, \label{en_2}
\end{equation}
for $j=1,...n-2$. 

And similarly, on the gridpoint directly upwind from the boundary, in each tube:
\begin{align}
\frac72 \lp \frac{T_{1,2}-T_{1,1}}{\Delta z}\rp -\frac{\Phi^2}{\rho_{1,2}^3}\lp\frac{\rho_{1,2}-\rho_{1,1}}{\Delta z}\rp+ \frac{\sigma'}{\Phi}\lp T_{1,2}-T_{2,2}\rp +g' =0,\label{en_1a}\\
\frac72 \lp \frac{T_{2,n}-T_{2,n-1}}{\Delta z}\rp -\frac{\Phi^2}{\rho_{2,n-1}^3}\lp\frac{\rho_{2,n}-\rho_{2,n-1}}{\Delta z}\rp+ \frac{\sigma'}{\Phi}\lp T_{1,n-1}-T_{2,n-1}\rp +g' =0. \label{en_2a}
\end{align}
Lastly, to complete the system, we have the boundary conditions:
\begin{align}
\frac72 T_{1,1}+\frac12 \frac{\Phi^2}{\rho_{1,1}^2} -\frac72 \Tint' = 0,\label{discbc1}\\
\frac{p_{1,1}}{\lp\rho_{1,1}\rp^{7/5}}-1 = 0,\label{discbc2}\\
\frac72 T_{2,n}+\frac12 \frac{\Phi^2}{\rho_{2,n}^2} - \frac72 = 0,\label{discbc3}\\
\frac{p_{2,n}}{\lp\rho_{2,n}\rp^\frac75} -1=0,\label{discbc4}\\
p_{1,n}-1=0,\label{discbc5}\\
p_{2,1}-\pint = 0. \label{discbc6}
\end{align}

Equations \eqref{mo_1a} through \eqref{discbc6}, together with the ideal gas law at every grid point, define a nonlinear system that we solve for the unknowns $\rho_{i,j}, T_{i,j}$ and $p_{i,j}$ for $i=1,2$ and $j = 1,...,n$, and also $\Phi$ and $\pint$.  Thus, after using the ideal gas law to remove $p_{i,j}$ from the list of unknowns we must solve for, the number of unknowns is
 $2(2n) + 2$. 
 After eliminating pressure from these equations, we are left with two equations connecting
 that grid point to its upwind neighborhood at each grid point. This count is applicable to the grid points at the entrance to each tube, provided that we consider the ``upwind neighborhood" of such a
 grid point to be its upstream reservoir, to which it is coupled by the compressible Bernoulli equation and also the adiabatic equation of state.  So far we have accounted for two equations per grid point, or $2(2n)$ overall.  The two additional equations that we need are provided by pressure continuity
 at the exit of each tube with the pressure in the air that is
 external to the tube at its exit. We solve these equations using Newton's method, as detailed below.
 
Let $\mathbf{X}$ be the vector of unknowns\footnote{We only store $\rho_{i,j}$ and $T_{i,j}$, since $p_{i,j} = \rho_{i,j}T_{i,j}$.},
 \begin{equation*}
     \mathbf{X} = \begin{bmatrix}
         \rho_{1,1} \ldots \rho_{1,n}& \rho_{2,1} \ldots \rho_{2,n}& T_{1,1} \ldots T_{1,n}& T_{2,1} \ldots T_{2,n}& \Phi &\pint
     \end{bmatrix},
 \end{equation*}
 $\mathbf{F}(\mathbf{X})$ be the function that evaluates the left hand sides of equations \eqref{mo_1} through \eqref{discbc6}, and $J$ be the Jacobian of $\mathbf{F}$ with respect to $\mathbf{X}$. At each step of Newton's method, we compute the $k$th guess for $\mathbf{X}^{(k)}$ by solving the following linear system, using Gaussian elimination:
 \begin{equation*}
     J\left(\mathbf{X}^{(k-1)}\right)\left(\mathbf{X}^{(k)}-\mathbf{X}^{(k)}\right) = - F\left(\mathbf{X}^{(k-1)}\right).
 \end{equation*}
We stop iterating when $\left\| F\left(\mathbf{X}^{(k)}\right)\right\|$ or $\left\|\mathbf{X}^{(k)}-\mathbf{X}^{(k)}\right\|$ are less than $10^{-12}$. 

As an initial guess for the Newton iteration, we use the following:

\begin{align}
\nonumber \rho_{1,j}=1/\Tint', \qquad \rho_{2,j} = 1,\\
\nonumber T_{1,j} = \Tint', \qquad T_{2,j} = 1,\\
\pint = 1,\qquad \Phi= 1/2.
\end{align}

%-----------------------------------------------------------------------------%

\section{Numerical Experiments} \label{numexp}
\begin{figure}[h!]
    \centering
    \begin{subfigure}{0.3\linewidth}
        \centering
        \includegraphics[width=\linewidth]{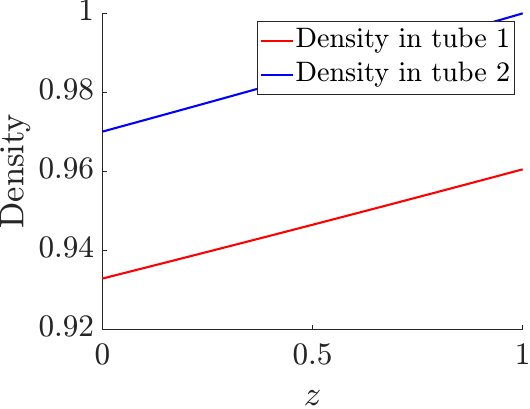}
    \end{subfigure}
        \begin{subfigure}{0.3\linewidth}
        \centering
        \includegraphics[width=\linewidth]{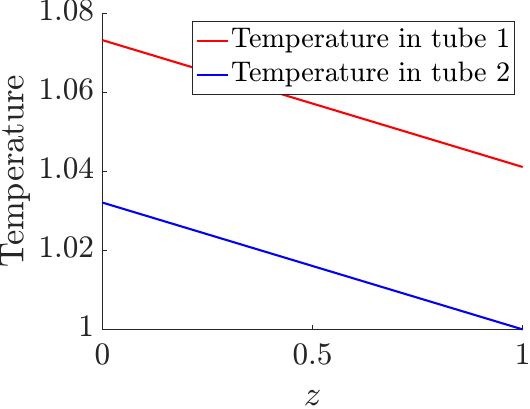}
    \end{subfigure}
        \begin{subfigure}{0.3\linewidth}
        \centering
        \includegraphics[width=\linewidth]{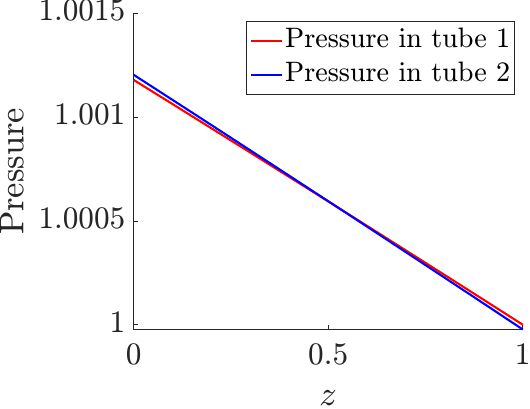}
    \end{subfigure}
    \caption{The dimensionless density, temperature, and pressure in both tubes when the model parameters are set to the physically realistic $g'=1.2471\times 10^{-4}$, $\Tint'=1.0732$, and $\sigma' = 0.0018$. }
    \label{fig:exsolution_std}
\end{figure}

\begin{figure}
    \centering
    \includegraphics[width=0.5\linewidth]{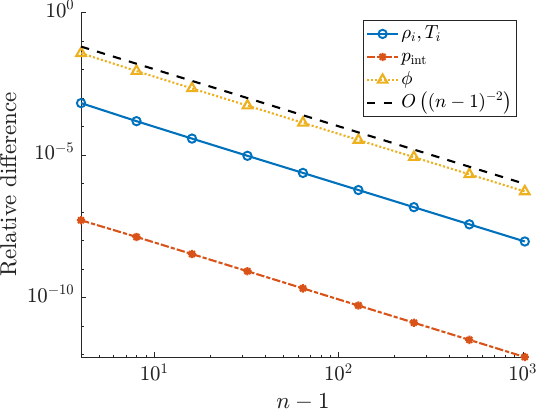}
   \caption{Relative maximum norm of the difference between computed solutions on successive grids during mesh refinement by a factor of 2, plotted on a log-log scale against $n-1$, the number of intervals on the coarser of the two grids, with the parameters set as above ($g'=1.2471\times 10^{-4}$, $\Tint'=1.0732$, and $\sigma' = 0.0018$.). The black dashed line indicates second order convergence. The curves for the relative maximum norm differences for the four variables $\rho_1, \rho_2, T_1$, and $T_2$ would be indistinguishable on this plot, and so the biggest relative difference only was plotted.}
    \label{fig:convstdy_std}
\end{figure}

 First, we take a look at the density, temperature, and pressure profiles at a physically reasonable choice of parameters (Figure \ref{fig:exsolution_std}). As shown in the figure, these variables look almost linear (though they are not) for this choice of parameters. It should be noted that a physically realistic $g'$ is small; we will come back to that in the next section. Figure \ref{fig:convstdy_std} shows the results of an empirical convergence study. At each step, the mesh width is refined by a factor of 2, so that the number of mesh points increases from $n$ to $2n-1$.  This ensures that each set of mesh points is a subset of the mesh points at the next level of refinement.  Then for any variable $v_i(z)$, we evaluate

 \begin{equation}
 \frac{\text{max}|v_i^n(z)-T_1^{2n-1}(z)|}{\text{max}|v_i^{2n-1}(z)|},
 \end{equation}
 where the maximum is taken over all values of $z$ that are common
to both meshes. The above relative difference is plotted in
Figure \ref{fig:convstdy_std} as a function of $n-1$ on a log-log scale.
%The log of the above relative difference is plotted in Figure \ref{fig:convstdy_std} as a function of $\log{n-1}$.  
For a second-order accurate method, the resulting plot should have slope -2, and a reference line with this slope is plotted for comparison.\footnote{For the two
special unknowns $\Phi$ and $\pint$ we use essentially the same formula, without the subscript $i$ and without the maximum over $z$: namely, $|v^n - v^{2n-1}| / |v^{2n-1}|$.}

Note in particular that second-order accuracy is achieved, despite the use of first-order accurate discretization at the mesh point
immediately downwind of the entry point of each tube.
 
\begin{figure}
    \centering
    \begin{subfigure}{0.3\linewidth}
        \centering
        \includegraphics[width=\linewidth]{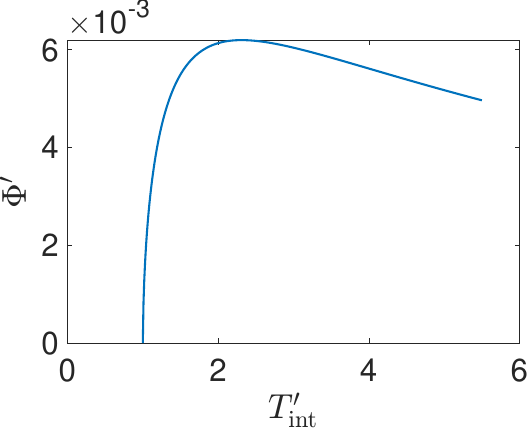}
    \end{subfigure}
        \begin{subfigure}{0.3\linewidth}
        \centering
        \includegraphics[width=\linewidth]{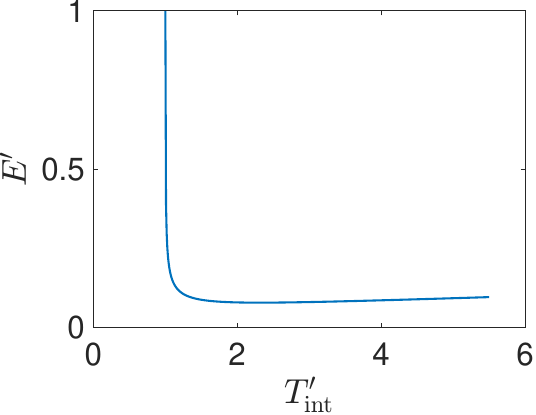}
    \end{subfigure}
        \begin{subfigure}{0.3\linewidth}
        \centering
        \includegraphics[width=\linewidth]{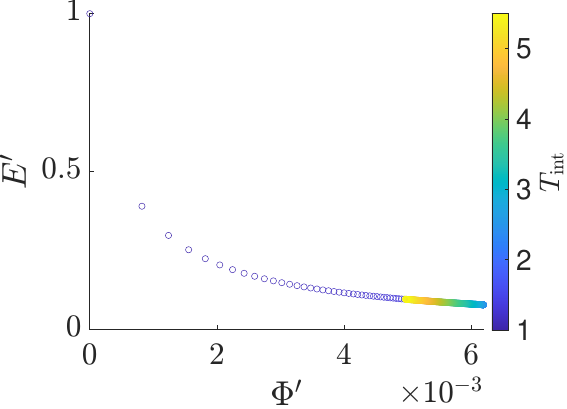}
    \end{subfigure}
    \caption{The mass flux and efficiency when varying $\Tint'$ and keeping the other parameters as above. We note that the efficiency is very high (and mass flux very low) when the temperature difference between the internal and external air is small. The efficiency decreases as $\Tint'$ varies through a range of physically realistic values and beyond, but starts to increase once again once the internal and external air are sufficiently far apart. Similarly, the mass flux increases then decreases as $\Tint'$ varies. Interestingly, the right-most plot, which tells us the trade-off between efficiency and mass flux, shows that the efficiency versus mass flux curve doubles back on itself as $\Tint'$ increases.}
    \label{fig:E_Phi_vsTintp}
\end{figure}

\begin{figure}
    \centering
    \begin{subfigure}{0.4\linewidth}
        \centering
        \includegraphics[width=\linewidth]{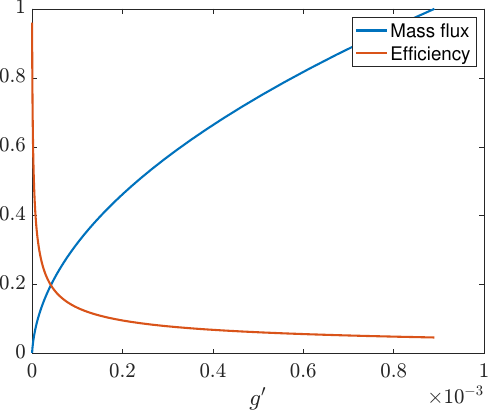}
    \end{subfigure}
        \begin{subfigure}{0.4\linewidth}
        \centering
        \includegraphics[width=\linewidth]{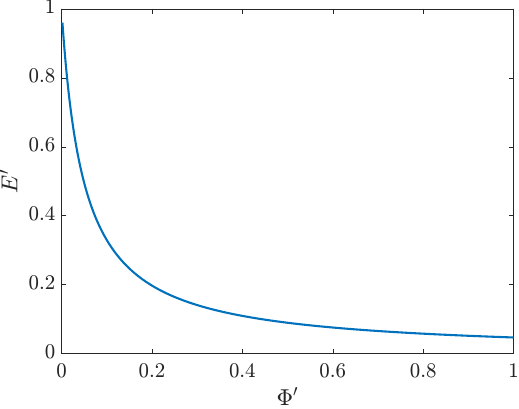}
    \end{subfigure}
    \caption{Mass flux and efficiency as the dimensionless parameter $g'$ changes.}
    \label{fig:E_Phi_vsgp}
\end{figure}

 Having looked at the physically realistic scenario, we then vary the parameters, starting with $\Tint'$. Figure \ref{fig:E_Phi_vsTintp} shows the mass flux and efficiency when varying $\Tint'$, and keeping the other parameters as above. Figure \ref{fig:E_Phi_vsgp} does the same, for varying $g'$.

\begin{figure}
    \centering
    \begin{subfigure}{0.4\linewidth}
        \centering
        \includegraphics[width=\linewidth]{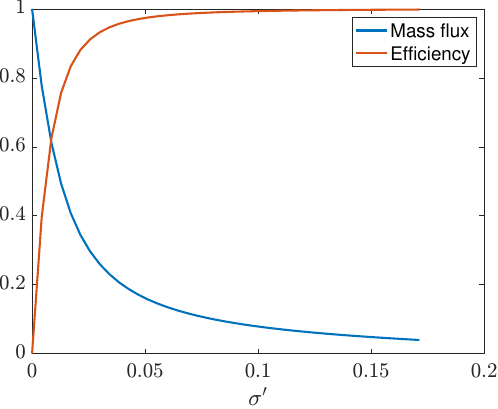}
    \end{subfigure}
        \begin{subfigure}{0.4\linewidth}
        \centering
        \includegraphics[width=\linewidth]{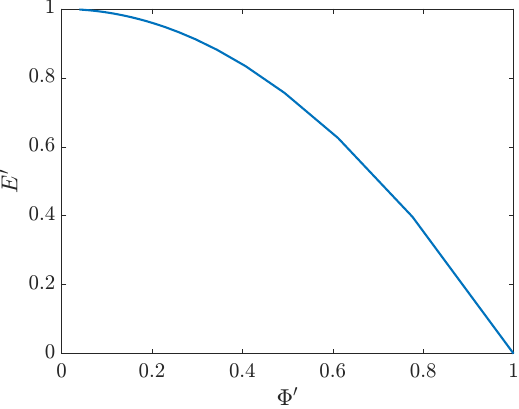}
    \end{subfigure}
    \caption{Relative mass flux ($\Phi_\text{rel}$) and efficiency ($E$) as functions of the dimensionless conductivity ($\sigma'$) of the partition, with $g'$ and $\Tint'$ held constant.  The parameter $\sigma'$ can be varied without changing the other dimensionless parameters by changing properties of the partition between the two pipes, such as its thickness or the thermal conductivity of the material out of which the partition is made. }
    \label{fig:E_Phi_vssigmap}
\end{figure}

 Figure \ref{fig:E_Phi_vssigmap} shows the tradeoff between efficiency and mass flux as we  vary the dimensionless conductivity $\sigma'$ of the partition.  Here, to get a measure of the mass flux that is comparable to efficiency, we normalize the mass flux by dividing by its maximum value: that is, the value of the mass flux when $\sigma = 0$ (so that the partition is thermally insulating and there is no heat exchange).  We use the notation $\Phi_{\text{rel}}$ for the mass flux that has been normalized in this way.  (Note that the relative mass flux is the same whether it is evaluated by normalizing the dimensionless mass flux or by normalizing
 the physical mass flux.)

From Figure \ref{fig:E_Phi_vssigmap}, we see clearly that the objective of maximizing the efficiency is in conflict with that of maximizing the relative mass flux.  The two curves in the left-hand panel of Figure \ref{fig:E_Phi_vssigmap} are the objective functions that we would like to maximize.  The minimum (worse-case) of the two objective functions is maximized where the two curves cross, so in that sense the crossing point determines the optimal value of $\sigma'$ for the design of a heat exchanger. It is encouraging that the efficiency and relative mass flux are about 0.6 at this point, which implies that we can achieve 60\% efficiency while simultaneously getting 60\% of the fresh air that would occur through the same pair of pipes with no heat exchange. 

One can argue, however, that a more sensible objective is to minimize the energy cost per unit mass of the fresh air that is being brought into the house.  By this criterion, one should make the efficiency as large as possible, even though the flow of fresh air is then very small, and then get as much fresh air as may be desired by using multiple copies of the highly efficient heat exchanger.  This approach is limited only by the number of heat exchangers that the roof can accommodate, or perhaps by the cost of the heat exchangers themselves. In any case, the curves of Figure \ref{fig:E_Phi_vssigmap} describe a relationship that needs to be considered in the design of a gravity-driven heat exchanger.

\begin{figure}
    \centering
    \begin{subfigure}{0.3\linewidth}
        \centering
        \includegraphics[width=\linewidth]{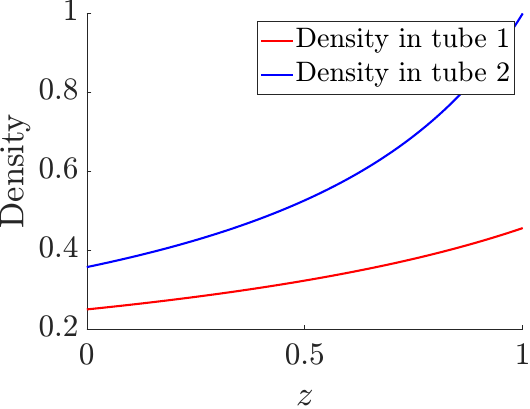}
    \end{subfigure}
        \begin{subfigure}{0.3\linewidth}
        \centering
        \includegraphics[width=\linewidth]{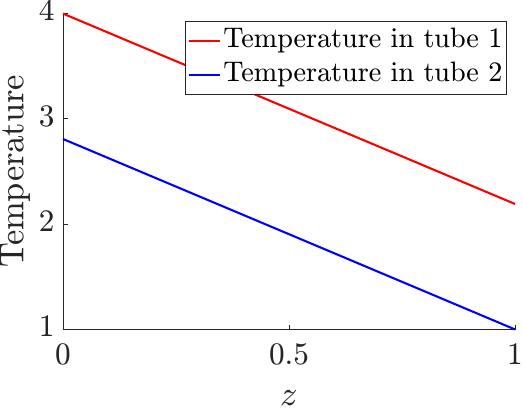}
    \end{subfigure}
        \begin{subfigure}{0.3\linewidth}
        \centering
        \includegraphics[width=\linewidth]{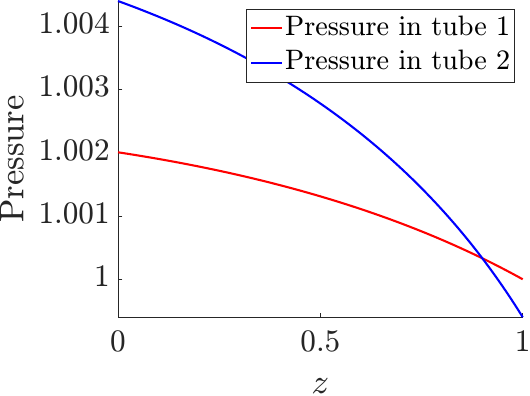}
    \end{subfigure}
    \caption{The dimensionless density, dimensionless temperature, and dimensionless pressure in both tubes, when the model parameters are set to the unrealistic values of $g'=0.0125$, $\Tint'=4$, and $\sigma' = 0.1836$. These values are chosen to give the method a stress test.}
    \label{fig:exsolution_high}
\end{figure}

\begin{figure}
    \centering
    \includegraphics[width=0.5\linewidth]{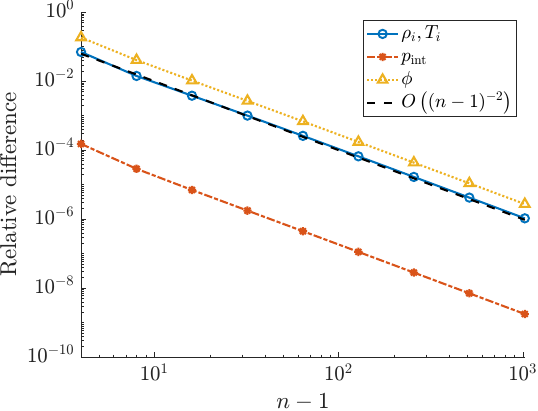}
    \caption{Relative maximum norm of the difference between computed solutions on successive grids during mesh refinement by a factor of 2, plotted on a log-log scale against $n-1$, the number of intervals on the coarser of the two grids, with the parameters set as above ($g'=0.0125$, $\Tint'=4$, and $\sigma' = 0.1836$.). The black dashed line indicates second order convergence. }
    \label{fig:convstdy_high}
\end{figure}

 Figure \ref{fig:exsolution_high} shows an example solution for a high value of $\Tint'$ and $\sigma'$, as a stress test for our method. It performs well at these values, as shown in Figure \ref{fig:convstdy_high}. In contrast to the previous results (Figure \ref{fig:exsolution_std}) shown for less extreme
 conditions, note in Figure \ref{fig:exsolution_high} the marked nonlinearity of the density
 and pressure distributions in tube 2.

\begin{figure}
    \centering
    \begin{subfigure}{0.4\linewidth}
        \centering
        \includegraphics[width=\linewidth]{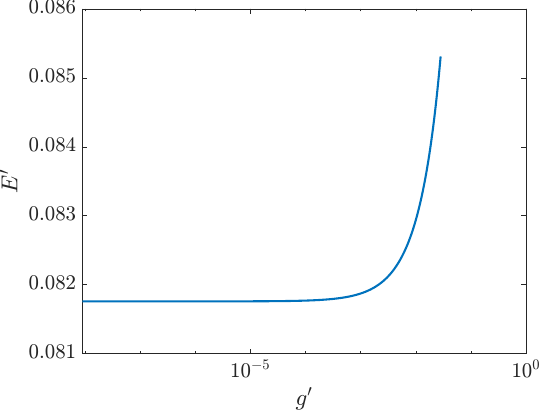}
    \end{subfigure}
        \begin{subfigure}{0.4\linewidth}
        \centering
        \includegraphics[width=\linewidth]{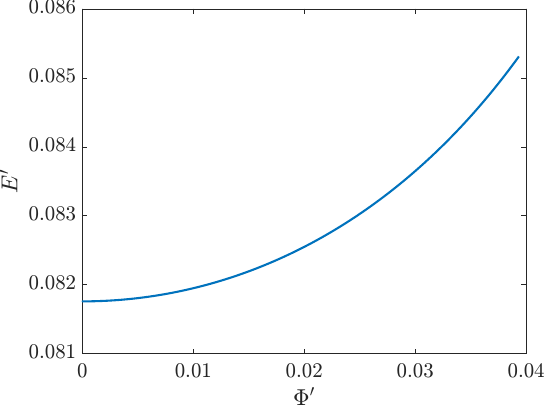}
    \end{subfigure}
    \caption{The efficiency plotted as a function of $g'$ (on the left) and the mass flux (on the right) as we vary both $g'$ and $\sigma'$, with $\sigma'$ proportional to $\sqrt{g'}$. Note that for $g'$ small, the efficiency barely changes. This prompts us to look at the asymptotic limit $g'\rightarrow 0$, with $\sigma'$ proportional to $\sqrt{g'}$.}
    \label{fig:E_Phi_varysigmapandgp}
    %plots in Phi_E_change_g.m
\end{figure}

 Finally, we note a relationship that emerged from playing with the
 parameters: when both $g'$ and $\sigma'$ are varied simultaneously, with $\sigma'\propto\sqrt{g'}$, for small $g'$ the efficiency is nearly constant (Figure \ref{fig:E_Phi_varysigmapandgp}). We explore this via asymptotic expansion in the next section.

\section{Asymptotic Predictions} \label{section5}

We expect the dimensionless parameter $g'$ to be small; indeed, even for a heat exchanger as tall as $O(10)$ meters and an outdoor temperature as cold as freezing, $273.15K$, $g'$ only gets as big as $O(0.001)$.

\subsection{Ansatz}	
We assume the following expansions in the small parameter $g'$:
\begin{align}
&\sigma' = \sqrt{g'}\sigma'^{(0)} \\
&T_i(z) = T_i^{(0)}(z)  + g'T_{i}^{(1)}(z)\\
&p_i(z) = p_i^{(0)}(z)  + g'p_{i}^{(1)}(z)\\
&\rho_i(z) = \rho_i^{(0)}(z)  + g'\rho_{i}^{(1)}(z)\\
&p_{\intt} = p_{\intt}^{(0)} + g'p_{\intt}^{(1)}\\
&\Phi = \sqrt{g'}\left( \Phi^{(0)}+g'\Phi^{(1)}\right). \end{align}
We want expressions for $T_i^{(0)}(z), p_i^{(0)}(z), \rho_i^{(0)}(z), p_{\intt}^{(0)},$ and $\Phi^{(0)}$.

\subsection{Conservation of Momentum, and Energy Dissipating Boundary Conditions}
Our conservation of momentum equation,
\begin{equation}
\frac{\partial}{\partial z} \left(p_i + \frac{\Phi^2}{\rho_i} \right) = -\rho_i g',
\end{equation}
gives us the following information:
\begin{align}
&\mathcal{O}(g'^0): \frac{\partial}{\partial z} \left(p_i^{(0)}\right) = 0 \\
&\mathcal{O}(g'): \frac{\partial}{\partial z} \left( p_i^{(1)} + \frac{(\Phi^{(0)})^2}{\rho_i^{(0)}}\right)=-\rho_i\nau.
\label{eq:momconsb}
\end{align}

So we have that, to lowest order (or, when $g'=0$), the pressure in each tube is constant. We can then use the following boundary conditions to get the constants:

\begin{equation}
p_1(1) = 1, \quad p_2(0)=p_{\intt}.
\label{eq:out_bc}
\end{equation}

This gives us that $p_1^{(0)}\equiv 1$ and $p_1^{(1)}(1)= 0$ in tube 1, as well as $p_2^{(0)}\equiv p_{\intt}^{(0)}$ and $p_2^{(1)}(0)=p_{\intt}^{(1)}$ in tube 2. 

Or, in other words:
\begin{align}
&p_1(z)=1+g'p_1^{(1)}(z), \quad\text{with }  p_1^{(1)}(1)= 0, \\
&p_2(z)=p_{\intt}^{(0)}+g'p_2^{(1)}(z), \quad\text{with } p_2^{(1)}(0)= p_{\intt}^{(1)}.
\end{align}

\subsection{Conservation of Energy, and Compressible Bernoulli Boundary Conditions}\label{sec:nrg_cons}

First, let's recall the conservation of energy equation:
\begin{equation}
\frac{\partial}{\partial z} \left(\frac{7}{2}T_i + \frac{1}{2}\frac{\Phi^2}{\rho_i^2}\right) = -\frac{\sigma'}{\Phi}(T_1-T_2)-g'.
\end{equation}
Substituting in to our expansion in $g'$, we find out that:
\begin{align}
&\mathcal{O}(g'^0): \frac{\partial}{\partial z} \left(\frac{7}{2}T_i^{(0}\right) = -\frac{\sigma'^{(0)}}{\Phi^{(0)}}T_D^{(0)} \\
&\mathcal{O}(g'): \frac{\partial}{\partial z} \left(\frac{7}{2}T_i^{(1)}+\frac{1}{2}\left(\frac{\Phi^{(0)}}{\rho_i^{(0)}}\right)^2 \right) = -\frac{\sigma'^{(0)}}{\Phi^{(0)}}\left(T_D^{(1)}-\frac{\Phi^{(1)}}{\Phi^{(0)}}T_D^{(0)}\right)-1,
\end{align}
where $T_D= T_1-T_2$.

We can use the boundary conditions
\begin{align}
&\frac{7}{2}T_1(0)+\frac{1}{2}\left(\frac{\Phi}{\rho_1(0)}\right)^2=\frac{7}{2}T_{\intt} \label{eq:T1_BC} \\
&\frac{7}{2}T_2(1)+\frac{1}{2}\left(\frac{\Phi}{\rho_2(1)}\right)^2=\frac{7}{2}\label{eq:T2_BC} ,
\end{align}

to get an expression for the temperatures at 0th order:
\begin{align}
&T_1^{(0)}(z)=-\frac{2}{7}\frac{\sigma'^{(0)}}{\Phi^{(0)}}\cdot T_D^{(0)}\cdot z+T_{\intt} \\
&T_2^{(0)}(z)=-\frac{2}{7}\frac{\sigma'^{(0)}}{\Phi^{(0)}}\cdot T_D^{(0)}\cdot (z-1)+1,
\end{align}
where $T_D=\frac{T_{\intt}'-1}{1+\frac{2}{7}\frac{\sigma'^{(0)}}{\Phi^{(0)}}}$.

\subsection{Remaining Boundary Conditions}

We have one remaining boundary condition in each tube, and have yet to determine $p_{\intt}^{(0)}$ and $\Phi^{(0)}$. Let's start with:
\begin{equation}
\frac{p_1(0)}{\rho_1(0)^{7/5}} = \frac{p_\intt}{\rho_{\intt}^{7/5}}
\label{eq:p1_bc}
\end{equation}

To 0th order, this gives us that $p_\intt\nau=1$.

Next, we have that
\begin{equation}
\frac{p_2(1)}{\rho_2(1)^{7/5}} = 1.
\label{eq:p2_bc}
\end{equation}

This does not give us anything new, so we need to find something else to determine $\Phi^{(0)}$.

\subsection{Using the Order $g$ part}
Let's revisit our boundary conditions, starting with the order $g'$ part of \eqref{eq:T1_BC} and \eqref{eq:T2_BC}:
\begin{align}
&\frac72 T_1^{(1)}(0)+\frac12 \lp \Phi\nau\Tint\rp^2 =0
\label{eq:t11_bc}\\
&\frac72 T_2^{(1)}(1)+\frac12 \lp \Phi\nau\rp^2 =0.
\label{eq:t21_bc}
\end{align}

Next, the order $g$ part of \eqref{eq:p1_bc} and \eqref{eq:p2_bc}:
\begin{align}
&\frac75\frac{T_1^{(1)}(0)}{\Tint}-\frac25p_1^{(1)}(0)=-\frac25\pint^{(1)}
\label{eq:p11_bc}\\
&\frac75 T_2^{(1)}(1)-\frac25p_2^{(1)}(1)=0.
\label{eq:p21_bc}
\end{align}

If we take twice \eqref{eq:t11_bc} and subtract $5\Tint\times$\eqref{eq:p11_bc}, we get the following:
\begin{equation}
\lp \Phi \nau \Tint' \rp ^2 + 2\Tint' p_1^{(1)}(0)=2 \Tint' \pint^{(1)}.
\end{equation}

Now we take twice \eqref{eq:t21_bc} minus $5\times$\eqref{eq:p21_bc}:
\begin{equation}
\lp \Phi \nau \rp ^2 + 2p_2^{(1)}(1)=0.
\end{equation}

Lastly, we look at the order $g$ term of the outflow boundary conditions:
\begin{equation}
p_1^{(1)}(1)=0, \quad p_2^{(1)}(0)=\pint^{(1)}.
\end{equation}

We have looked at all the boundary conditions, and now we use the order $g'$ term of the momentum conservation equation \eqref{eq:momconsb}. It tells us that
\begin{align}
&\frac{\partial}{\partial z} \left( p_i^{(1)} + (\Phi^{(0)})^2\frac{T_i\nau}{p_i^{(0)}}\right)=-\rho_i\nau \\
\implies &\frac{\partial}{\partial z} \lp p_i^{(1)}\rp = -(\Phi^{(0)})^2 \frac{\partial}{\partial z} \lp  T_i \nau \rp -\rho_i\nau,
\end{align}
since $p_i\nau = 1$, and $\Phi\nau$ is constant. Taking the integral thus gives
\begin{equation}
p_i^{(1)}(1)-p_i^{(1)}(0)= -\lp\Phi\nau\rp^2\lp T_i\nau(1)-T_i\nau(0)\rp-\int_0^1\rho_i\nau dz.
\end{equation}

We know the pressures at the inflow and outflow boundaries to 0th order, as well as the temperatures. If we use this information for tube 1, we can write:
\begin{align}
&0-\lp\pint\one - \frac{\Tint'}{2}(\Phi\nau)^2\rp = -(\Phi\nau)^2\lp-\frac{2}{7}\frac{(\sigma')\nau}{\Phi\nau}T_D\nau\rp-\int_0^1 \frac{1}{T_1\nau (z)} dz \\
\implies &\pint\one = (\Phi\nau)^2\lp\frac{\Tint'}{2}-\frac27\frac{(\sigma')\nau}{\Phi\nau}T_D\nau\rp + \int_0^1 \frac{1}{-\frac{2}{7}\frac{\sigma'^{(0)}}{\Phi^{(0)}} T_D^{(0)}z+T_{\intt}} dz.
\label{eq:alpha}
\end{align}

Similarly, in tube 2:
\begin{align}
&-\frac12(\Phi\nau)^2 - \pint\one = -(\Phi\nau)^2\lp-\frac{2}{7}\frac{(\sigma')\nau}{\Phi\nau}T_D\nau\rp-\int_0^1 \frac{1}{T_2\nau (z)} dz \\
\implies &\pint\one = (\Phi\nau)^2\lp-\frac12-\frac27\frac{(\sigma')\nau}{\Phi\nau}T_D\nau\rp + \int_0^1 \frac{1}{-\frac{2}{7}\frac{(\sigma')^{(0)}}{\Phi^{(0)}} T_D^{(0)} (z-1)+1} dz. 
\label{eq:beta}
\end{align}

If we set the right hand sides of \eqref{eq:alpha} and \eqref{eq:beta} equal, we get an equation to solve for $\Phi\nau$:
\begin{gather}
(\Phi\nau)^2\lp\frac{\Tint'}{2}-\frac27\frac{(\sigma')\nau}{\Phi\nau}T_D\nau + \frac12+\frac27\frac{(\sigma')\nau}{\Phi\nau}T_D\nau\rp \qquad \qquad
\nonumber \\ =  \int_0^1 \frac{1}{-\frac{2}{7}\frac{(\sigma')^{(0)}}{\Phi^{(0)}} T_D^{(0)} (z-1)+1} dz -  \int_0^1 \frac{1}{-\frac{2}{7}\frac{(\sigma')^{(0)}}{\Phi^{(0)}} T_D^{(0)}z+T_{\intt}} dz.
\end{gather}
Thus,
\begin{align}
\nonumber (\Phi\nau)^2\lp\frac{\Tint'+1}{2}\rp &= \left. \frac{\log\lp-\frac27\frac{(\sigma')^{(0)}}{\Phi^{(0)}} T_D^{(0)} (z-1)+1\rp}{-\frac27\frac{(\sigma')^{(0)}}{{\Phi^{(0)}}}T_D^{(0)}}\right|_{z=0}^1 - \left.\frac{\log\lp-\frac27\frac{(\sigma')^{(0)}}{\Phi^{(0)}} T_D^{(0)}z +\Tint'\rp}{-\frac27\frac{(\sigma')^{(0)}}{{\Phi^{(0)}}}T_D^{(0)}}\right|_{z=0}^1 \\
&= \lp\frac72\frac{\Phi\nau}{(\sigma')\nau}\frac{1}{T_D\nau}\rp\log\lp(1+\frac27\frac{(\sigma')\nau}{\Phi\nau}T_D\nau)(1-\frac27\frac{(\sigma')\nau}{\Phi\nau}\frac{T_D\nau}{\Tint'})\rp.
\end{align}

We can then use a root solver, such as MATLAB fzero, to get $\Phi\nau$, and this gives us all the 0th order terms. 

\begin{figure}
    \centering
    \includegraphics[width=0.5\linewidth]{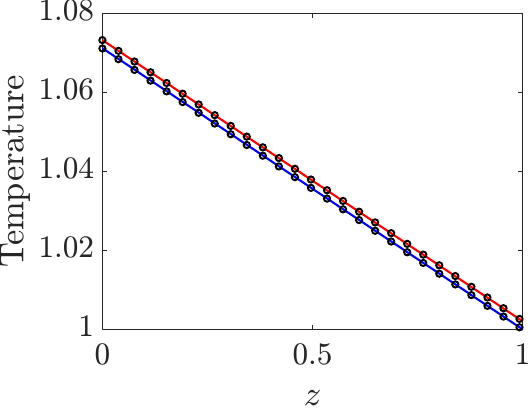}
    \caption{The numerically computed temperature in tube 1 and 2, respectively in red and blue, with the asymptotic solution overlaid in black circles. The numerical solution was calculated using a realistic value for $g'$, and we see excellent agreement between the asymptotic solution and this numerically computed solution.}
    \label{fig:comparenumasy}
\end{figure}

Once we have those, we can compare the numerical results to the asymptotic ones. For the purpose of this comparison, we choose the thermal conductivity to be 3 $Wm^{-1}K^{-1}$ (for context, aluminum has a thermal conductivity of 237 $Wm^{-1} K^{-1}$, aluminum oxide of 30 $Wm^{-1}K^{-1}$, and air of about 0.02 $W m^{-1}K^{-1}$ at standard temperature and pressure; so we choose something between these values, somewhat arbitrarily pending a more rigorous study in an upcoming article). For this choice of thermal conductivity and corresponding value of $\sigma'$, the numerically computed $\Phi$ converges to $\Phi\nau$. Figure \ref{fig:comparenumasy} shows the close correspondence between the asymptotic temperature and the numerically computed temperature in both tubes, for a physically realistic choice of $g'$.
%$\Phi\nau$ agrees with 

\section{Summary and Conclusions}
In this paper, we introduce a complete model of an air-to-air, buoyancy driven heat exchanger, based on one-dimensional gas dynamics. We consider only the steady state equations, and devise a second-order accurate numerical method to explore the consequences of these equations. An asymptotic method is also introduced, based on a very small parameter, the ratio of the height of the heat exchanger to the thickness of the atmosphere in the adiabatic atmospheric model. The asymptotic and numerical results are in excellent agreement.

Making use of these methods, we are able to evaluate the efficiency of the heat exchanger. We define efficiency by considering the energy required to heat each unit mass of fresh air that is brought into the house. In the absence of heat exchange, this is simply the energy cost of raising the temperature of the unit mass of air from the external to the internal temperature. The relative reduction in this power requirement in the presence of heat exchange is our definition of efficiency (see equation \eqref{effdef}).

We find that there is a tradeoff in the design of a buoyancy driven heat exchanger between minimizing heat loss and maximizing fresh air. As the thermal conductivity of the partition between the inflow and outflow increases, heat loss is reduced, but air flow is reduced as well. If we normalize the flow of fresh air by the largest amount that could flow through the heat exchanger, that is, by the amount of fresh air that would flow through it in the absence of heat exchange, then the optimal tradeoff is achieved when the normalized air flow is equal to the efficiency of the heat exchanger. At this point, with reasonable choices of parameters, the efficiency of the heat exchanger is approximately 60\%.

\section{Appendix}
In this appendix, we argue that an energy-dissipating boundary condition at the outflow of each tube is necessary by looking at the consequences of conserving energy throughout the entire system. We consider here only the simplest special case, in which there is no heat transferred across the partition; in other words, we consider for each tube a column of air in adiabatic equilibrium, acted upon only by gravity.

In steady flow, the following quantity is then constant along streamlines:
\begin{equation}
    \frac72 \frac{RT_i(\xi)}{m} + \frac12 u_i^2(\xi) + g\xi,\label{conservedquantity}
\end{equation}
where $\xi$ is the height at any point in the flow.

Since there is a streamline connecting the reservoirs inside and outside of the building, the quantity \eqref{conservedquantity} must be equal in the two reservoirs. Since the velocity is zero by hypothesis in each reservoir, this in turn means that
\begin{equation}
    \frac72 \frac{R}{m}\lp\Text-\Tint\rp + gH = 0,\label{superspecialcase}
\end{equation}
or, in other words, the system is overdetermined and there can only be a solution in a very special case: the case in which relation \eqref{superspecialcase} between $\Tint$, $\Text$, and $H$ holds true. Note that this is a very small temperature difference; for example, if $H = 10$m, then $\Tint-\Text=0.0973$K.

\nocite{*}
\bibliographystyle{plain}
\bibliography{biblio}

\end{document}